\documentclass[a4 paper, 11pt]{smfart}

\input macros

\title[{\sc Minimum essentiel et degrŽs d'obstruction des translat{\'e}s de sous-tores}]
{\uppercase{Minimum essentiel et degrŽs d'obstruction\\ des translat{\'e}s de sous-tores}}

\author{Patrice {\sc Philippon}}

\address{Institut de Math{\'e}matiques de Jussieu (U.M.R. 7586),
Projet G{\'e}om{\'e}trie et Dynamique.
Case~7012,
2 place Jussieu,
75251 Paris Cedex 05,
France.}
\email{pph@math.jussieu.fr}
\urladdr{http://www.math.jussieu.fr/\~{}pph/}

\author{Mart{\'\i}n {\sc Sombra}}

\address{Universitat de Barcelona,
Departament d'{\`A}lgebra i Geometria.
Gran Via 585,
08007 Barcelona, Espagne.}
\email{sombra@ub.edu}
\urladdr{http://atlas.mat.ub.es/personals/sombra/}
\thanks{ M. Sombra a \'et\'e financ\'e par le {Programme Ram\'on y Cajal} du Ministerio de Educaci\'on y Ciencia, Espagne.}



\begin{abstract}
Nous Žtudions les degrŽs d'obstruction des translatŽs de sous-tores des tores multiplicatifs et nous montrons comment ils interviennent dans les minorations 
du  minimum essentiel de ces mmes variŽtŽs. En particulier, nous combinons nos calculs avec des rŽsultats de A.~Schinzel et de F.~Amoroso - R.~Dvornicich, que nous Žtendons ainsi aux translatŽs de sous-tores dŽfinis sur des corps CM, montrant que leurs minimums essentiels ont un comportement gŽomŽtrique de mme nature que ceux des variŽtŽs qui ne sont pas translatŽes de sous-tores. \'Egalement, nous raffinons et gŽnŽralisons ˆ des translatŽs de sous-tores 
une minoration due ˆ D.~Bertrand en liaison avec les relations de dŽpendance multiplicative des points algŽbriques. 

Nous montrons enfin que, dans le cas fonctionnel la minoration conjecturale du minimum essentiel des sous-variŽtŽs des tores multiplicatifs est vraie et que sa dŽmonstration est ŽlŽmentaire.
\end{abstract}

\begin{altabstract}
{\bf  Essential minimum and obstruction degrees of translates of sub-tori.}  We study the obstruction degrees of translates of sub-tori of multiplicative tori and we show how they are connected to lower bounds for the essential minimum of these varieties. In particular, we combine our computations with results of A.~Schinzel and of  F.~Amoroso - R.~Dvornicich, that we thus extend to translates of sub-tori defined over  CM fields. This shows that the essential minimum of these varieties have a geometric behavior similar to that of varieties which are not translates of sub-tori. We also refine and extend to translates of sub-tori a lower bound obtained by D.~Bertrand in 
connexion with relations of multiplicative dependence of algebraic points. 

Finaly, we prove that in the functional case the conjectural lower bound 
for the essential minimum of sub-varieties of multiplicative tori is true 
and its proof is elementary.
\end{altabstract}
\begin{document}

\maketitle


\setlength{\baselineskip}{13pt}

\vspace{-8mm}
\section{Introduction et rŽsultats}

Du point de vue de la gŽomŽtrie diophantienne, les translatŽs de 
sous-tores du tore multiplicatif $\G_m^N:=(\Qbar^\times)^{N}$ se 
trouvent ˆ la croisŽe des problmes de Lehmer et de Bogomolov g{\'e}n{\'e}ralis{\'e}s. 
En effet, on sait que lorsqu'elles ne sont pas de torsion, les minorations pour la hauteur normalisŽe de ces variŽtŽs sont de nature arithmŽtique, dŽpendant essen\-tiel\-le\-ment du corps de dŽfinition de la variŽtŽ. Au contraire, pour les sous-variŽtŽs 
qui ne sont pas des translatŽs de sous-tores, on dispose de minorations ne dŽpendant que de leur gŽomŽtrie, {\it voir}~les 
conjectures~\ref{sinnou} et~\ref{sinnoubogo} ci-dessous, {\it voir} aussi~\cite{DP99,AD01,AD03,Dav02}.

\bigskip

Le degrŽ d'une variŽtŽ n'est pas le bon invariant pour l'Žtude des minorations de 
la hauteur, il faut considŽrer plut™t son degrŽ d'obstruction. Soient $X \subset Y$ des ensembles algŽbriques Žquidimensionnels de $\P^N$ dŽfinis sur un corps algŽbriquement clos $\k$ et $K$ un sous-corps de $\k$. Le {\it degrŽ d'obstruction de $X$ dans $Y$ relatif ˆ $K$ }
est 
$$
\omega_K(X;Y):= \min \left\{\deg(W) \, : \ W\in {\rm Div}_K(Y), W\supset X \right\} 
\enspace,
$$
\cad le plus petit degr{\'e} d'un diviseur $W$ de $Y$ d{\'e}fini sur $K$ et contenant $X$, 
vu comme un cycle de $\P^N$. 
On convient de poser $\omega_K(X;Y):=\infty$ lorsque qu'aucun diviseur de $Y$ dŽfini sur $K$ ne contient $X$. 

La plus petite hauteur d'un sous-ensemble de points algŽbriques Zariski dense  
dans une variŽtŽ est quantifiŽe par la notion de minimum essentiel.  
Soit $X \subset \P^N$  une vari{\'e}t{\'e} quasi-projective quelconque dŽfinie sur $\Qbar$, et pour $ \theta\ge 0$ posons $X(\theta)$ l'ensemble des points algŽbriques de $X$ de hauteur de Weil (ou {\it hauteur normalisŽe}) $\hnorm$ majorŽe par $\theta$. Le {\it minimum essentiel de $\hnorm$ sur $X$} est 
$$
\wh{\mu}^\ess(X) := \inf \{\theta : X(\theta) \mbox{ est Zariski dense}\}
\enspace. 
$$
Aussi, on dŽsigne par $\wh{\mu}^\abs(X)$ le {\it minimum absolu de $\hnorm$ 
sur $X$}, c'est l'infimum des hauteurs de points algŽbriques de $X$. 

\smallskip

On fixera la compactification \'equivariante  de $\G_m^N$ donn{\'e}e par l'inclusion 
standard $\iota:~\G_m^N\hookrightarrow \P^N, (t_1,\dots,t_N)\mapsto (1:t_1:\cdots:t_N)$. Ceci permet de transporter les notions de hauteur, minimum essentiel, etc., aux sous-objets de $\G_m^N$ en considŽrant leur adhŽrence de Zariski dans $\P^N$.

Un {\em sous-tore} $T\subset {\bf G}_m^N$ est un sous-groupe algŽbrique isomorphe ˆ un ${\bf G}_m^n$. 
Rappelons que tous les sous-groupes algŽbriques connexes de $\G_m^N$ sont de 
cette forme, et qu'en gŽnŽral les sous-groupes algŽbriques de ${\bf G}_m^N$ 
sont les produits $F\cdot T$ d'un groupe fini $F$ par un sous-tore $T$. Nous notons $\mug^\infty$ le sous-groupe de $\G_m$ des racines de l'unit{\'e}.
Une {\it variŽtŽ de torsion} $U\subset \G_m^N$ est le translatŽ d'un sous-tore par un point de torsion, \cad un point ˆ coordonnŽes dans $\mug^\infty$.
Pour une variŽtŽ $X\subset \G_m^N$ on pose $U_X$ la plus petite variŽtŽ de torsion
la contenant.

\medskip 
Notre motivation de dŽpart est la minoration conjecturale suivante pour la hauteur des points Zariski dense dans une variŽtŽ. Elle est proposŽe par F.~Amoroso et S.~David, en s'inspirant de divers rŽsultats sur les problme de Lehmer et de Bogomolov gŽnŽralisŽs; 
c'est une rŽŽcriture de~\cite[Conj.~1.10]{AD03} ({\it voir} \S~\ref{degobst}):

\begin{conj} \label{sinnou}
Il existe un rŽel $c(N)>0$ tel que pour toute variŽtŽ $X\subset \G_m^N(\Qbar)$ et $U_X^\Q$ la plus petite rŽunion de variŽtŽs de torsion contenant $X$ qui soit dŽfinie sur $\Q$, on ait \footnote{Si $\dim(X)=\dim(U^\Q_X)$ on a $\omega_\Q(X,U^\Q_X)=\infty$ et $\frac{\deg(U^\Q_X)}{\omega_\Q(X;U^\Q_X)}=0$.}
$$
\munorm^\ess(X) \geq c(N) \frac{\deg(U^\Q_X)}{\omega_\Q(X;U^\Q_X)}\enspace.
$$
\end{conj}
\smallskip

Cette conjecture est Žquivalente, ˆ la valeur de la constante $c(N)$ prs, ˆ son analogue sur la cl™ture abelienne $\Q^\ab$ de $\Q$ (proposition~\ref{sinnouabelien})~:
$$
\munorm^\ess(X) \geq c'(N) \frac{\deg(U_X)}{\omega_{\Q^\ab}(X;U_X)}\enspace.
\enspace$$

Elle est surtout int{\'e}ressante pour les translatŽs de sous-tores par des points d'ordre infini, car pour les autres variŽtŽs on s'attend ˆ une minoration ˆ caractre gŽo\-mŽ\-trique, ne dŽpendant pas du corps de dŽfinition de la variŽtŽ (\og problme de Bogomolov effectif\fg): 

\begin{conj}[\cite{AD03}, Conj.~1.2] \label{sinnoubogo}
Il existe un rŽel $c(N)>0$ tel que pour toute variŽtŽ $X\subset \G_m^N(\Qbar)$ et $V_X$ le plus petit translatŽ de sous-tore contenant $X$, on ait
$$
\munorm^\ess(X) \geq c(N) \frac{\deg(V_X)}{\omega_\Qbar(X;V_X)}\enspace.
$$
\end{conj}

Notons toutefois que cette minoration, si elle est de nature purement gŽomŽtrique, 
n'en est pas toujours meilleure que celle de la conjecture~\ref{sinnou}, {\it voir}
remarque~\ref{diffomegas}.

\medskip 
Pour un sous-corps  $K \subset \Qbar$ on pose
$$
m(K) := 
\inf\{ \hnorm(\alpha) \, : \ \alpha\in K^\times \setminus {\mug}^\infty\}\enspace.
$$
\smallskip
Dans ce texte nous prŽsentons la minoration suivante pour la hauteur des points 
d'un translatŽ de sous-tore par un point d'ordre infini:

\begin{thm}\label{laurier}

Soit $X\subset \G_m^N$ un translatŽ de sous-tore, dŽfini sur un corps de nombres $K$, 
alors 
$$
\munorm^\abs(X) \geq c_1(N) \cdot m(K) \cdot 
 \frac{\deg(U_X)}{\omega_{\overline{\Q}}(X;U_X)} \enspace, 
\quad \mbox{ avec}\quad 
\declare{ctedslaurier}(N) \ge 
N^{-3/2} 2^{-N}\enspace.
$$
\end{thm}

Ce thŽorme concerne le minimum {absolu} de $\hnorm$ sur $X$. 
Or, pour les translat{\'e}s de sous-tores les minimums absolu et 
essentiel co{\"\i}ncident et donc cette minoration peut de faon Žquivalente 
se reformuler en termes de $\munorm^\ess$. 
Notons aussi qu'un translatŽ de sous-tore $X$ est dŽfini sur $K$ si et 
seulement si c'est le translatŽ par un point dŽfini sur $K$, \cad si et seulement s'il contient un point $K$-rationnel, {\it voir} \S~\ref{degobstetres}. De plus, la variŽtŽ de torsion minimale $U_X$ est automatiquement dŽfinie sur $K$ elle aussi, {\it voir} lemme~\ref{U_XestsurK}.

\smallskip



\medskip
Le probl{\`e}me de Lehmer classique consiste {\`a} montrer l'existence d'une constante $c>0$ telle que $m(K)\geq c [K:{\bf Q}]^{-1}$ pour tout corps de nombres $K$. Le meilleur r{\'e}sultat connu dans cette direction reste celui de E.~Dobrowolski~\cite{Dob79}
$$
m(K) 
\geq \frac{c}{[K:\Q]} \left(\frac{\log\log(3[K:\Q])}{\log(2[K:\Q])}\right)^3
\enspace.
$$
De fait, Amoroso et David ont dŽmontrŽ les conjectures~\ref{sinnou} et \ref{sinnoubogo} 
ˆ des facteurs logarithmiques prs, par une extension non triviale de l'approche de  Dobrowolski~\cite{AD01,AD03}.

Par ailleurs, on sait par des travaux de A.~Schinzel~\cite{Sch73} que lorsque $K$ est un corps totalement rŽel ou CM (c'est-ˆ-dire une extension imaginaire quadratique d'un corps totalement rŽel) et $\alpha\in K^\times$ satisfait $|\alpha|\not=1$, alors
\begin{equation}\label{schinzel}
\hnorm(\alpha) \geq \frac{1}{2}\log\left(\frac{1+\sqrt{5}}{2}\right)\enspace,
\end{equation}
avec ŽgalitŽ si et seulement si $\alpha=\frac{1\pm\sqrt{5}}{2}$. 
D'aprs~\cite{AN05}, la condition $|\alpha|\not=1$ ne peut tre ŽliminŽe, notons cependant que lorsque $\alpha$ est entier mais pas une racine de l'unitŽ, elle est satisfaite par au moins un conjuguŽ de $\alpha$ et lorsque $K$ est totalement rŽel elle est satisfaite par tout $\alpha\in K^\times\setminus\mug^\infty$. Amoroso et R.~Dvornicich ont montrŽ~\cite{ADv00} qu'ˆ l'instar du cas totalement rŽel, dans une extension abelienne $K$ de ${\bf Q}$ on a une minoration uniforme pour la hauteur des nombres algŽbriques non nuls qui ne sont pas des racines de l'unitŽ~:
\begin{equation}\label{amodvo}
m(K) \geq \frac{\log(5)}{12}\enspace.
\end{equation}

En combinant le thŽorme~\ref{laurier} avec les inŽgalitŽs~(\ref{schinzel}) et~(\ref{amodvo}) on obtient la gŽnŽralisation suivante  
(ˆ la constante prs) des rŽsultats de~\cite{Sch73} et~\cite{ADv00}~:

\begin{cor}\label{corCM}
Soit $X \subset \G_m^N$ un translatŽ de sous-tore dŽfini sur un corps  $K$ totalement rŽel ou une extension abelienne de ${\bf Q}$, alors 
$$
\munorm^\abs(X) \geq \declare{ctedscorCM}(N) \frac{\deg(U_X)}{\omega_{\overline{\bf Q}}(X;U_X)} \enspace, 
\quad \mbox{ avec}\quad\use{ctedscorCM}(N):=
\frac{\log(5)}{12} \use{ctedslaurier}(N) \enspace.
$$
Plus gŽnŽralement, si $K$ est un corps CM notons $G_X$ la plus petite variŽtŽ de torsion 
contenant $X\cdot\overline{X}$ ($\overline{X}$ dŽsigne ici l'image de $X$ sous l'action de la conjugaison complexe), alors
$$\munorm^\abs(X) \geq \use{ctedscorCM}(N)  \frac{\deg(G_X)}{\omega_{\overline{\bf Q}}(X\cdot\overline{X};G_X)}
\enspace.$$
\end{cor}

On montre que, dans la situation de ce  corollaire et pour un corps CM, la variŽtŽ de torsion $G_X$ est en fait un sous-tore (lemme~\ref{translatecasCM}). La minoration du corollaire~\ref{corCM} pour $K$ un corps totalement rŽel se compare favorablement ˆ celle de la conjecture~\ref{sinnou}, gr‰ce au lemme~\ref{deQaK}. Dans le cas d'une extension abelienne de $\Q$, elle est par ailleurs consŽquence de la conjecture~\ref{sinnou} suivant un argument de David, {\it voir} proposition~\ref{sinnouabelien}.

\bigskip
La preuve de ces rŽsultats s'appuie sur une Žtude des degrŽs d'obstruction des translatŽs de sous-tores dans les tores. Remarquons que si $X \subset Y$ sont des variŽtŽs 
dŽfinies sur un sous-corps $K \subset \Qbar$, il n'existe pas nŽcessairement de diviseur dŽfini sur $K$ rŽalisant $\omega_\Qbar(X;Y)$ et on n'a donc pas en gŽnŽral 
$\omega_K(X;Y)=\omega_\Qbar(X;Y)$ comme le montrent des exemples au \S~\ref{degobst}. C'est toutefois le cas lorsque $Y$ est un translatŽ de sous-tore, {\it voir} corollaire~\ref{cosette}.

On dŽmontre que lorsque $X$ et $Y$ sont des translatŽs de sous-tores de dimensions respectives $n$ et $p$, le degrŽ d'obstruction relatif ˆ $\Qbar$ se rŽalise par une Žquation binomiale, 
et qu'il est comparable au premier minimum d'un certain rŽseau (proposition~\ref{corbeau}). Comme consŽquence, on peut le majorer gr‰ce au premier thŽorme de 
Minkowski par (corollaire~\ref{omega})
$$
\omega_\Qbar(X;Y) \leq 4^{N}{N}^{3/2} \deg(X)^{\frac{1}{p-n}} \deg(Y)^{1-\frac{1}{p-n}} \enspace.
$$
L'inŽgalitŽ plus prŽcise du lemme~\ref{vautour} jointe au thŽorme~\ref{laurier}
entra{\^\i}ne~:
\begin{cor}\label{bertrand}
Soit $X\subset \G_m^N$ un translatŽ de sous-tore par un point d'ordre infini et 
dŽfini sur un corps de nombres $K$, alors 
$$
\wh{\mu}^\abs(X) \geq 2^{-N}N^{-2} \cdot m(K) 
\cdot \bigg(\frac{\deg(U_X)}{\deg(X)}\bigg)^{1/{{\rm codim}_{U_X}(X)}} \enspace.
$$
\end{cor}

Cette minoration peut alternativement se regarder comme une majoration pour le 
degrŽ de la variŽtŽ de torsion minimale $U_X$. Si $X$ est rŽduite ˆ un point 
$\alpha\in(K^\times)^N\setminus (\mug^\infty)^N$ 
$$
\hnorm(\alpha) \geq 
2^{-N}N^{-2} \cdot m(K) \cdot 
\deg(U_{\alpha})^{1 /\dim(U_\alpha)}\enspace.
$$
On retrouve ainsi la minoration (*) de~\cite{Ber95} pour le cas des tores dŽployŽs, avec de plus une dŽpendance prŽcise dans le corps $K$ (comparer avec~\cite[Cor.~1]{Ber95}). Ce type de minoration a des applications au calcul du module des relations de dŽpendance multiplicative des nombres, {\it voir} aussi ˆ ce propos~\cite{Ber97}.

\bigskip
Finalement, on Žtudie l'analogue fonctionnel de la conjecture~\ref{sinnou}. Soit $\k$ un corps al\-gŽ\-bri\-que\-ment clos, on s'intŽressera ˆ la hauteur des points dans une variŽtŽ $X\subset \P^N(\ov{\k(t)})$. Le r™le des racines de l'unitŽ est jouŽ par les ŽlŽments de $\k$, et celui des variŽtŽs de torsion par les variŽtŽs dŽfinies sur $\k$. 
Notons $U_X$ la plus petite variŽtŽ dŽfinie sur $\k$ et contenant $X$. 
On montre que si $X$ elle mme n'est pas dŽfinie sur $\k$, alors elle est de codimension 1 dans $U_X$, et que son minimum essentiel est minorŽ par (proposition~\ref{sinnoufonctionnel})
$$
\mu^{\ess}(X) \geq \frac{1}{d+1} \frac{\deg(U_X)}{\omega_{{\bf k}(t)}(X;U_X)}\enspace.
$$
La dŽmonstration se base sur la thŽorie de l'intersection multiprojective ŽlŽmentaire.

\bigskip
Le texte est organisŽ de la faon suivante~: aux~\S~\ref{degobst} et~\S~\ref{degobstetres} on Žtudie le degrŽ d'obstruction dans les sous-tores et on Žtablit sa relation avec le premier minimum d'un certain rŽseau. Nous montrons le thŽorme~\ref{laurier} et ses corollaires au~\S~\ref{demonstrations}, les dŽmonstrations n'utilisent que les arguments de nature gŽomŽtrique dŽveloppŽs aux paragraphes prŽcŽdents. Au~\S~\ref{lehmerfonct} nous Žtudions l'analogue fonctionnel de la conjecture~\ref{sinnou}. 

\bigskip
\noindent{\bf Remerciements.\hspace{1.5mm}---}\hspace{1mm} 
Nous remercions Sinnou David pour son aide au cours de la gŽnse de ce travail, ainsi que Francesco Amoroso pour ses remarques sur une premire version de ce texte.

\section{DegrŽs d'obstruction dans les tores}\label{degobst}

Le but de ce paragraphe est de dŽcrire en termes combinatoires les degrŽs des diviseurs des translatŽs de sous-tores. Pour cette Žtude on peut se ramener ˆ supposer, quitte ˆ translater la situation, que le translatŽ de sous-tore en question est en fait un sous-tore.

\medskip

\subsection{La situation gŽomŽtrique}\label{situgeom}
Soit $\k$ un corps algŽbriquement clos. On pose $\G_m^N(\k):=(\k^\times)^N$ le groupe multiplicatif et $\P^N(\k)$ l'espace projectif sur $\k$; on omettra la mention au corps lorsque celui-ci sera clair dans le contexte. Une {\it variŽtŽ} sera une sous-variŽtŽ rŽduite et irrŽductible de $\P^N(\k)$ ou de $\G_m^N(\k)$. Pour un sous-corps $K$ de $\k$, une $K$-variŽtŽ sera un sous-ensemble algŽbrique de $\P^N(\k)$ ou de $\G_m^N(\k)$, dŽfini sur $K$, rŽduit et $K$-irrŽductible. Pour une variŽtŽ $X$ on dŽsigne par $I(X)$ son idŽal de dŽfinition dans $ \k[x_0,\dots,x_N]$ ou dans $\k[x_1^{\pm 1},\dots, x_N^{\pm1}]$, suivant que $X\subset \P^N$ ou $X\subset \G_m^N$. 
\medskip

Commenons par dŽcrire les sous-tores ˆ l'aide de paramŽtrisations. Soit ${\cB}=(b_1,\dots,b_N)\in (\Z^p)^{N}$; le sous-tore associŽ $T_\cB \subset \G_m^N(\k)$ est par dŽfinition l'image de l'application monomiale
$$\varphi_\cB: \G_m^p(\k)\to \G_m^{N}(\k) \quad, \quad t \mapsto (t^{b_1} , \dots,
t^{b_N})\enspace.$$
On supposera dans la suite que le ${\bf Z}$-module $L_{\cB}:={\bf Z}b_1+\dots+{\bf Z}b_N$ co•ncide avec $\Z^p$. Ceci Žquivaut ˆ ce que $\varphi_\cB$ soit un isomorphisme entre 
$\G_m^p$ et $T_\cB$, et entra"ne donc $\dim(T_\cB)=p$. Pour $\beta=(\beta_1,\dots,\beta_N)\in\G_m^N$ on pose $\beta T_\cB$ le translatŽ de $T_\cB$ par le point $\beta$. 
Dans les notations de~\cite{PS04}, $T_\cB$ (resp. $\beta T_\cB$) est l'ouvert principal de la vari{\'e}t{\'e} torique $X_{\cB'} \subset {\bf P}^N$ (resp. $X_{\cB',\beta'}$) 
pour $ \cB':=({\mathbf 0},b_1,\dots,b_N)\in (\Z^p)^{N+1}$ et $\beta':=(1,\beta_1,\dots,\beta_N)\in (\k^\times)^{N+1}$. 

Soit $Q_\cB\subset \R^p$ l'enveloppe convexe des points 
${\mathbf 0},b_1,\dots,b_N$ et introduisons aussi l'application linŽaire 
$$
M_\cB:\R^N\to \R^p \quad , \quad  \lambda\mapsto \lambda_1 b_1+\cdots +\lambda_N b_N \enspace.
$$ 
Pour un polyn™me de Laurent 
$f\in\k[x_1^{\pm1},\dots,x_N^{\pm1}]$ on note 
$\Newton(f) \subset \R^N$ son {\it polytope de Newton}, enveloppe convexe des exposants des mon™mes apparaissant effectivement dans l'Žcriture de $f$.
Soit $\varphi_\cB^*(f):= f\circ \varphi_\cB 
\in\k[t_1^{\pm1},\dots,t_p^{\pm1}]$, alors
$ \Newton(\varphi^*_\cB(f))= M_\cB(\Newton(f))$. 

\medskip
Le {\it volume mixte} (ou {\it multi-volume}) $\MV(Q_1, \dots, Q_p)$ d'une famille
d'ensembles convexes $Q_1, \dots, Q_p \subset \R^p$ est dŽfini par la formule  
\begin{equation*} 
\MV(Q_1, \dots, Q_p) := 
\sum_{j=1}^p (-1)^{p - j} 
\sum_{1 \le i_1 < \cdots < i_j \le p} 
\Vol_n(Q_{i_1} + \cdots + Q_{i_j}) \enspace. 
\end{equation*} 
Cette notion g{\'e}n{\'e}ralise le volume d'un ensemble convexe, car on a $\MV(Q, \dots, Q) = n! \, \Vol_n(Q)$ o $\Vol_n$ dŽsigne le volume euclidien. Le volume mixte est sym{\'e}trique, lin{\'e}aire en chaque variable $Q_i$ par rapport {\`a} la somme de Minkowski et monotone par rapport ˆ l'inclusion. On renvoie {\`a}~\cite{Ewa96} pour les propri{\'e}t{\'e}s de base de cette notion. 

\medskip 
Le rŽsultat suivant donne une description combinatoire 
du degrŽ du diviseur dŽcoupŽ par $f$ sur $T_\cB$: 

\begin{prop}\label{degdiv}
Soit $\cB \in (\Z^p)^{N+1}$ tel que $L_\cB=\Z^p$ et $f\in \k[t_1^{\pm 1}, \dots, t_N^{\pm 1}] \setminus I(T_\cB)$, alors
$$\deg(\div(f)\cdot T_\cB) = 
\MV(M_\cB(\Newton(f)), \underbrace{Q_\cB,\dots,Q_\cB}_{p-1\ fois})\enspace.
$$
Plus gŽnŽralement, soit $Y=\bigcup_{i=1}^D \beta_i T_\cB$ une rŽunion de $D$ translatŽs de $T_\cB$, alors
$$
\deg(\div(f)\cdot Y) = D\cdot \MV(M_\cB(\Newton(f)), {Q_\cB,\dots,Q_\cB})\enspace.
$$
\end{prop}

\begin{demo}
Posons $W:=\div(f)\cdot T_\cB \in \Div(T_\cB)$, vu comme un cycle de dimension 
$p-1$ de ${\bf P}^N$. On a
$
\deg(W) = \Card\Big(W \cap Z(\ell_1,\dots,\ell_{p-1})\Big)
$ pour $\ell_i\in\k[x_1,\dots,x_N]$ des formes linŽaires gŽnŽriques et donc 
$$ 
\deg(W) = \Card \Big(\varphi_\cB^*(f) \cap  \varphi_\cB^*(\ell_1) \cap  
\cdots \cap \varphi_\cB^*(\ell_{p-1})\Big) 
$$
car $\varphi_\cB:\G_m^p \to T_\cB$ est une bijection. 
Le thŽorme de Bernstein-Koushnirenko~\cite{Ber75} entra"ne alors la majoration 
\begin{equation}\label{formuleBK}
\deg(W) \le \MV\left(\Newton({\varphi_\cB^*(f)}), \Newton(\varphi_\cB^*(\ell_{1})) , \dots, \Newton(\varphi_\cB^*(\ell_{p-1}))\right) = \MV\left(M_\cB(\Newton(f)), {Q_\cB,\dots,Q_\cB}\right)
\end{equation}
car $\Newton(\varphi_\cB^*(\ell_{i}))= Q_\cB$.  

\smallskip
Pour $\tau\in \R^p$ et $g\in \k[t_1^{\pm 1},\dots, t_p^{\pm 1}]
$ notons 
$\init_\tau(g)\in \k[t_1^{\pm 1},\dots, t_p^{\pm 1}]$ la partie initiale de $g$ relative au poids $\tau$, \cad, la partie de $g$ de poids maximal dans la direction de $\tau$.  
Bernstein dŽmontre aussi que la majoration~(\ref{formuleBK}) est une ŽgalitŽ si et seulement si 
pour tout  $\tau\in \R^p$ le systme 
$$
\init_\tau (\varphi_\cB^*(f))=0 \enspace, \quad \init_\tau (\varphi_\cB^*(\ell_1))=0 \enspace, \enspace \dots \enspace, \enspace \init_\tau (\varphi_\cB^*(\ell_{p-1}))=0
$$
n'a pas de solution dans $\G_m^p$. Or, $\init_\tau (\varphi_\cB^*(f))$ est toujours un polyn™me de Laurent non nul puisque $f\ne 0$, et donc cette condition est satisfaite par des formes linŽaires $\ell_i$ gŽnŽriques, ce qui dŽmontre le rŽsultat.

\smallskip 
Pour vŽrifier la dernire ŽgalitŽ de la proposition, par additivitŽ du degrŽ il suffit de vŽrifier 
$$\deg(\div(f)\cdot T_\cB) = \deg(\div(f)\cdot (\beta T_{\cB}))$$
pour tout $\beta\in\G_m^N(\k)$. On a $\div(f)\cdot (\beta T_{\cB})=
\beta\cdot (\div((\beta^{-1})^*(f))\cdot T_\cB)$. 
Or le degrŽ est invariant par translation et $\varphi_\cB^*\circ(\beta^{-1})^*(f)$ et $\varphi_\cB^*(f)$ ont mme support, l'identitŽ cherchŽe rŽsulte donc de la formule prŽcŽdemment Žtablie.
\end{demo}

\begin{cor} \label{cosette}
Soit $K \subset \k$ un sous-corps et $X\subset T_\cB$ une 
$K$-variŽtŽ, alors dans les notations de la proposition~\ref{degdiv} 
$$
\omega_K(X; T_\cB) = \min \Big\{ \MV\left(M_\cB(\Newton(f)) , Q_\cB,\dots,Q_\cB\right) \, : \ f\in I(X) \setminus I(T_\cB)\Big\}
\enspace.$$
En particulier $\omega_K(X;T_\cB) = \omega_{\k}(X;T_\cB)$.
\end{cor}

\begin{demo}
L'anneau $\k[T_\cB]$ est factoriel car isomorphe {\`a} $\k[t_1^{{\pm1}},\dots,t_p^{{\pm1}}]$. Il est donc principal et tout diviseur de $T_\cB$ contenant $X$ est dŽfini par une seule Žquation  
$f\in I(X) \setminus I(T_\cB
) $. 

On Žcrit $f$ comme combinaison linŽaire de mon™mes 
indŽpendants modulo $I(T_\cB)$. Comme $X$ est une $K$-variŽtŽ, la condition Žnonant que $f$ s'annule sur $X$ s'Žcrit au travers d'un systme linŽaire dŽfini sur $K$, on peut donc en trouver une solution non triviale $g\in K[x_1^{\pm1},\dots, x_N^{\pm1}]$ qui satisfait $\Supp(g) \subset \Supp(f)$ (donc 
$M_\cB(\Newton(g))\subset M_\cB(\Newton(f))$) 
et qui dŽfinit un $K$-diviseur de $T_\cB$. Il suit de la proposition~\ref{degdiv} et de la monotonie du volume mixte par rapport ˆ l'inclusion en chaque variable~\cite[Thm.~4.12]{Ewa96} que le degrŽ du $K$-diviseur 
$\div(g)\cdot T_\cB$ 
est infŽrieur ou Žgal ˆ celui de $\div(f)\cdot T_\cB$.  
Finalement, en considŽrant un diviseur rŽalisant $\omega_{\k}(X; T_\cB)$ on conclut $\omega_{\k}(X;T_\cB) = \omega_K(X;T_\cB)$, puis avec la proposition~\ref{degdiv} on obtient l'expression de $\omega_K(X;T_\cB)$ mentionnŽe.
\end{demo}

Lorsque $Y$ n'est pas un translatŽ de sous-tore, tout en restant dŽfini sur $K$, on n'a plus nŽcessairement l'ŽgalitŽ $\omega_K(X;Y) = \omega_{\k}(X;Y)$, mme lorsque $X$ est elle-mme un translatŽ de sous-tore dŽfini sur $K$.
En effet, considŽrons $X:=\{(1:1:1:1)\}$ dans la rŽunion $Y$ de quatre plans de ${\bf P}^3$ dŽfinie par l'Žquation, rationnelle sur ${\bf Q}$, suivante
$$\prod_{i,j\in\{0,1\}}\left((w-x)+(-1)^i\sqrt{2}(w-y)+(-1)^j\sqrt{3}(w-z)\right) = 0
\enspace.$$
Une droite rationnelle sur ${\bf Q}$ passant par le point $X$ admet une paramŽtrisation de la forme
$$
\P^1 \to\P^3 \enspace, \quad (t_0:t_1)\mapsto (t_0:at_1+t_0:bt_1+t_0:ct_1+t_0)
$$
avec $a,b,c\in{\bf Q}$ non tous nuls. On en dŽduit qu'une telle droite ne peut tre contenue dans $Y$, car $1$, $\sqrt{2}$ et $\sqrt{3}$ sont linŽairement indŽpendants sur ${\bf Q}$, et donc $\omega_{\bf Q}(X;Y) \geq 2$ (en fait $=2$), tandis que $\omega_{\k}(X;Y) = 1$ car $Y$ contient une infinitŽ de droites (dŽfinies sur $\k$) passant par $X$.

On peut aussi produire un exemple o $Y$ est irrŽductible: soit $Y\subset \P^3$ la surface, dŽfinie sur ${\bf Q}$ et irrŽductible sur $\k$, d'Žquation $wz-x^2+2y^2 = 0$ 
et $X$ un des points $(1:0:0:0)$ ou $(0:0:0:1)$. Cette surface contient deux droites dŽfinies sur ${\bf Q}(\sqrt{2})$ passant par chacun des points fixŽs, mais aucune droite dŽfinie sur ${\bf Q}$ passant par un de ces mmes points.

\medskip 

\subsection{ParticularitŽs arithmŽtiques}\label{partarith}
On se place maintenant sur $\k=\Qbar$, soit $X \subset \G_m^N(\Qbar)$ un ensemble algŽbrique Žquidimensionnel et $K$ un corps de nombres. Posons $U_X^K$ la plus petite rŽunion de variŽtŽs de torsion contenant $X$ qui soit dŽfinie sur $K$. Notons qu'une rŽunion de variŽtŽs de torsion est toujours dŽfinie sur $\Q^\ab$, la cl™ture abelienne de $\Q$, et donc $U_X^K$ est dŽfinie sur $K\cap \Q^\ab$, la sous-extension abelienne maximale de $K$.

\begin{lem}\label{U_XestsurK}
Soit $X \subset \G_m^N(\Qbar)$ une variŽtŽ dŽfinie sur un corps de nombres $K$, alors $U_X=U_X^\Qbar$, la plus petite variŽtŽ de torsion contenant $X$, est dŽfinie sur $K\cap\Q^\ab$, elle co•ncide avec $U_X^K$. 
\end{lem}

\begin{demo}
Ceci rŽsulte de ce que la classe des rŽunions de variŽtŽs de torsion est fermŽe par 
intersections: soit $U_X$ la plus petite rŽunion de variŽtŽs de torsion contenant $X$,
alors $\sigma(U_X) =U_X$ pour tout $\sigma \in \Gal(\Qbar/K)$, donc $U_X$ est dŽfinie sur $K$ et par suite sur $K\cap\Q^\ab$. 
\end{demo}

Comme consŽquence de ceci et de l'ŽgalitŽ $\omega_K(X;Y) = \omega_K\left(\bigcup_{\sigma\in{\rm Gal}(\overline{\bf Q}/K)}\sigma(X);Y\right)$, la conjecture~\ref{sinnou} peut se reformuler de la faon suivante, plus proche de la formulation originale dans~\cite[Conj.~1.10]{AD03}: soit $X\subset \G_m^N$ une $\Q$-variŽtŽ et posons $U_X$ la plus petite rŽunion de variŽtŽs de torsion contenant $X$, alors il existe un r\'eel $c(N)>0$ indŽpendant de $X$ tel que
$$
\munorm^\ess(X) \geq c(N) \frac{\deg(U_X)}{\omega_\Q(X;U_X)}\enspace.
$$

\begin{lem} \label{deQaK}
Soit $X \subset \G_m^N$ une variŽtŽ et $K_1$ le corps de dŽfinition de $U_X$, alors
\begin{eqnarray*}
\deg(U_X^\Q) &= &[K_1:\Q]\deg(U_X)\enspace,\\[0mm]
\omega_\Q(X;U_{X}^\Q) &= & [K_1:\Q] \omega_{K_1}(X;U_X)\enspace.  
\end{eqnarray*}
En particulier, si $X$ est dŽfinie sur un corps de nombres $K$ on a
$$
\frac{1}{[K:K_1]}\cdot\frac{\deg(U_X)}{\omega_\Qbar(X;U_X)} \le \frac{\deg(U_X^\Q)}{\omega_{\Q}(X;U_X^\Q)}  = \frac{\deg(U_X)}{\omega_{K_1}(X;U_X)} \le \frac{\deg(U_X)}{\omega_{\Qbar}(X;U_X)}
\enspace.$$
\end{lem}

\begin{demo}
On a $U_X^\Q= \bigcup_{\sigma} \sigma(U_X)$ o $\sigma$ parcourt les $\Q$-plongements de $K_1$ dans $\Qbar$, d'o la premire ŽgalitŽ. Soit $E \subset U_X$ un $K_1$-diviseur rŽalisant $\omega_{K_1}(X;U_X)$, alors $\bigcup_{\sigma} \sigma(E)$ est un $\Q$-diviseur de $U_X^\Q$ de degrŽ $\le [K_1:\Q] \omega_{K_1}(X;U_X)$ et contenant $X$, d'o $\omega_{\Q}(X;U_X^\Q)\le[K_1:\Q] \omega_{K_1}(X;U_X)$. Si maintenant $E^\Q$ est un $\Q$-diviseur de $U_X^\Q$ rŽalisant $\omega_\Q(X;U_X^\Q)$ alors $E^\Q\cap U_X$ est un $K_1$-diviseur de $U_X$ contenant $X$ et donc de degrŽ $\ge\omega_{K_1}(X;U_X)$. Mais $\bigcup_{\sigma} \sigma(E^\Q\cap U_X)$ est contenu dans $E^\Q$ et de degrŽ $\ge [K_1:\Q]\omega_{K_1}(X;U_X)$, on a donc $\omega_\Q(X;U_{X}^\Q) \ge [K_1:\Q] \omega_{K_1}(X;U_X)$ et la seconde ŽgalitŽ.
\smallskip

L'ŽgalitŽ centrale de la dernire formule rŽsulte de ce qui prŽcde et les inŽgalitŽs extrmes de l'encadrement
$$\omega_{\Qbar}(X;U_X)\le\omega_{K_1}(X;U_X)\le[K:K_1]\omega_{\Qbar}(X;U_X)
\enspace.$$
L'inŽgalitŽ de droite vient du fait qu'il existe un diviseur rŽalisant $\omega_{\Qbar}(X;U_X)$ dŽfini sur $K$, d'aprs le corolaire~\ref{cosette} et le lemme~\ref{U_XestsurK}.
\end{demo}

\begin{question}
Sous rŽserve d'une rŽponse positive au problme de Lehmer abelien le terme $m(K)\cdot\frac{\deg(U_X)}{\omega_\Qbar(X;U_X)}$, donnant la minoration dans le thŽorme~\ref{laurier}, est compris dans l'intervalle dŽterminŽ par les termes extrmes de la dernire formule du lemme~\ref{deQaK}. Comment se situe cette quantitŽ par rapport ˆ $\frac{\deg(U_X^\Q)}{\omega_{\Q}(X;U_X^\Q)}$, qui est Žgalement dans ce mme intervalle?
\end{question}


Le problme de Lehmer abelien consiste ˆ montrer l'existence d'un rŽel $c>0$ tel que $m(K^\ab)\geq c[K^\ab: \Q^\ab]^{-1}$, o $K^\ab$ dŽsigne la cl™ture abelienne de $K$. Dans cette direction, Amoroso et U.~Zannier~\cite{AZ00} ont montrŽ en direction du problme de Lehmer abelien que dans le thŽorme de Dobrowolski on peut essentiellement remplacer le degr{\'e} de $K$ sur $\Q$ par le degrŽ  $[K^\ab: \Q^\ab]$. RŽcemment Amoroso et E.~Delsinne~\cite{ADe06} on amŽliorŽ cette minoration en 
$$
m(K^\ab)\geq \frac{c}{[K^\ab: \Q^\ab]} \frac{\log \log(3 [K^\ab: \Q^\ab])^{3}}{\log(2 [K^\ab: \Q^\ab])^{4}} \enspace.
$$ 
Suivant un argument que nous a expliquŽ David, on peut montrer que la conjecture~\ref{sinnou} entra"ne une rŽponse positive au problme de Lehmer abelien. Nous gŽnŽralisons cet argument pour montrer que la conjecture~\ref{sinnou} est Žquivalente, ˆ la constante $c(N)$ prs, ˆ son analogue sur la cl™ture abelienne de $\Q$~: 

\begin{prop} \label{sinnouabelien}
Soit $X \subset \G_m^N$ une variŽtŽ dŽfinie sur $\Qbar$, sous la conjecture~\ref{sinnou} on a
\begin{equation}\label{conjabel}
\wh\mu^\ess(X) \geq c(N+1)\frac{\deg(U_X)}{\omega_{\Q^\ab}(X;U_X)}\enspace.
\end{equation}
RŽciproquement, la minoration (\ref{conjabel}) implique $\wh\mu^\ess(X) \geq c(N+1)\frac{\deg(U_X^\Q)}{\omega_{\Q}(X;U_X^\Q)}$.
\end{prop}

\begin{demo}
Soit $K$ un corps de nombres sur lequel $X$ et $U_X$ sont dŽfinis, et soit 
$\xi\in \mug^\infty$ telle que la variŽtŽ de torsion $U_X$ soit dŽfinie sur l'extension cyclotomique $\Q(\xi)$. Soit $L$ le compositum de $K$ et $\Q(\xi)$, considŽrons $\wt X := X\times\{\xi\} \subset\G_m^{N+1}$ qui est une variŽtŽ dŽfinie sur $L$. Posons la variŽtŽ de torsion
$$
Y:= \bigcup_{\sigma\in \gal(\Q(\xi)/\Q)} 
\sigma\Big(U_X \times \{\xi\}\Big) \enspace,
$$
qui, par contruction, est dŽfinie sur $\Q$, on a $ U_{\wt X}^\Q\subset Y$. En outre, il est clair que $U_X \times \{\xi\} \subset U_{\wt X}^\Q$ et en prenant la cl™ture galoisienne on obtient $ Y=U_{\wt X}^\Q$. On a donc
\begin{eqnarray*}
\deg(U_{\wt X}^\Q) 
&=& [\Q(\xi):\Q] \deg(U_X) \enspace,\\[2mm]
\omega_\Q(\wt X, U_{\wt X}^\Q) 
&= & [\Q(\xi):\Q] \omega_{\Q(\xi)}( \wt X, U_X\times\{\xi\}) = 
[\Q(\xi):\Q] \omega_{\Q(\xi)}( X, U_X) \enspace.  
\end{eqnarray*}
Et donc, en appliquant la conjecture~\ref{sinnou} ˆ $\wt X$, il vient
$$
\munorm^\ess(X)= \munorm^\ess(\wt X) \ge c(N+1) 
\frac{\deg(U_{\wt X}^\Q)}{\omega_\Q(\wt X, U_{\wt X}^\Q) }
= c(N+1) \frac{\deg(U_X)}{\omega_{\Q(\xi)}( X, U_X)} \enspace.  
$$
On conclut en prenant $\xi$ d'ordre suffisamment grand pour que 
$\omega_{\Q(\xi)}( X, U_X)=\omega_{\Q^\ab}( X, U_X)$; c'est possible car 
$\Q^\ab=\Q(\mug^\infty)$ (thŽorme de Kronecker-Weber). 
\smallskip

Pour la rŽciproque on remarque que d'aprs le lemme~\ref{deQaK}, si $K_1\subset\Q^\ab$ dŽsigne le corps de dŽfinition de $U_X$, on a 
$$\frac{\deg(U_X^\Q)}{\omega_{\Q}(X;U_X^\Q)} = \frac{\deg(U_X)}{\omega_{K_1}(X;U_X)} \leq \frac{\deg(U_X)}{\omega_{\Q^\ab}(X;U_X)}\enspace.$$
\end{demo}

\section{DegrŽs d'obstruction et rŽseaux}\label{degobstetres}

Dans ce paragraphe on Žtudie le degrŽ d'obstruction d'un translatŽ de sous-tore dans un autre sous-tore, relativement ˆ un corps algŽbriquement clos $\k$. Soit $\cB=(b_1,\dots,b_N)\in (\Z^p)^N$ tel que $L_\cB=\Z^p$ et, comme precŽdemment, posons $T_\cB \subset \G_m^N$ le sous-tore associŽ, image de l'application 
$$
\varphi_\cB: \G_m^p\to \G_m^N \quad, \quad s \mapsto (s^{b_1}, \dots, s^{b_N}) \enspace,
$$
et $\beta T_\cB$ son translat{\'e} par un point $\beta \in \G_m^N$. Le noyau $\ker(M_\cB)\subset \Z^N$ est un sous-module de rang $N-n$, {\it saturŽ}, \cad tel que le quotient $\Z^N/\ker(M_\cB)$ soit sans torsion. Les idŽaux de dŽfinition de $T_\cB$ et de $\beta T_\cB$ dans $\k[x_1^{\pm1},\dots,x_N^{\pm1}]$ s'Žcrivent~\cite{ES96}
\begin{equation}\label{moineau}
I(T_\cB) = \Big( x^\lambda - 1 \, : \ \lambda \in \ker(M_\cB) \Big) 
\quad \mbox{ et } \quad  
I(\beta T_\cB) 
= \Big( x^\lambda - \beta^\lambda \, : \ \lambda \in \ker(M_\cB) \Big)  
\enspace.\end{equation}
Ainsi, le sous-tore $T_\cB$  ne dŽpend que de $\ker(M_\cB)$. Plus gŽnŽralement, le translatŽ $\beta T_\cB$ ne dŽpend que de $\ker(M_\cB)$ et de l'homomorphisme $\ker(M_\cB)\to \k^\times, \lambda\mapsto \beta^\lambda$. 

\medskip
RŽciproquement, soit  $\Gamma \subset \Z^N$ un sous-module saturŽ muni d'un caractre $\rho$, \cad un homomorphisme $\rho:\Gamma\to \k^\times$. Ces donnŽes dŽfinissent un idŽal binomial
$$
I(\Gamma,\rho):= 
\Big( x^\lambda - \rho(\lambda) \, : \ \lambda \in \Gamma \Big) 
\enspace \subset \k[x_1^{\pm1},\dots,x_N^{\pm1}] \enspace.
$$
On peut vŽrifier que tout sous-module saturŽ de $\Z^N$ se rŽalise comme $\ker(M_\cB)$ pour un certain $ \cB \in (\Z^p)^N$ du type envisagŽ ci-dessus, et que tout caractre d'un tel module se rŽalise par un point $\beta \in \G_m^N$ (on peut mme choisir $\beta \in \rho(\Gamma)^N$). Les correspondances 
$$
(\Gamma,\rho) \mapsto I(\Gamma,\rho) \quad \mbox{ et } \quad I \mapsto Z(I)
$$
sont des bijections entre les sous-modules saturŽs $\Gamma$ de $\Z^N$ munis d'un caractre  $\rho$, les idŽaux premiers binomiaux de $\k[x_1^{\pm1},\dots,x_N^{\pm1}]$ ne contenant aucune des variables $x_j$, et les translatŽs de sous-tores de $\G_m^N$~\cite[Cor.~2.6]{ES96}.

On Žcrira aussi $T(\Gamma,\rho)$ pour le translatŽ de sous-tore associŽ au couple
$\Gamma, \rho$. Le caractre trivial $\rho=1$ correspond ˆ un sous-tore et ceux pour lesquels $\rho(\Gamma)\subset \mug^\infty$ aux variŽtŽs de torsion. RŽciproquement, pour un translatŽ de sous-tore $X \subset \G_m^N$ on dŽsigne par $\Gamma_X $  et $\rho_X$ le sous-module et le caractre  associŽs. La reprŽsentation par sous-modules saturŽs et caractres correspond donc aux Žquations des translatŽs de sous-tores. Dans ce paragraphe on privilŽgiera cette reprŽsentation, sauf dans la proposition~\ref{carpincho}. 

\medskip
Le translatŽ $T(\Gamma, \rho)$ est dŽfini sur un sous-corps $K$ si et seulement si 
$\rho(\Gamma) \subset K^\times$. Ceci Žquivaut ˆ ce qu'il existe $\beta  \in (K^\times)^N$ tel que $\rho(\lambda)= \beta^\lambda$ pour tout $\lambda\in\Gamma$, \cad 
ˆ ce que $T(\Gamma,\rho) $  contienne un point $K$-rationnel. En effet, comme $\Gamma$ est saturŽ l'homomorphisme $\rho$ s'Žtend en $\rho:\Z^N\rightarrow K^\times$, posant alors $\beta_1,\dots,\beta_N\in K^\times$ les images par $\rho$ des ŽlŽments de la base canonique de $\Z^N$ on vŽrifie que $\rho$ s'Žcrit bien $\lambda\mapsto\beta^\lambda$.

Soient $X,Y \subset \G_m^N$ des translatŽs de sous-tores, la condition $X\subset Y$ Žquivaut ˆ 
$$\Gamma_X \supset \Gamma_Y
\quad \mbox{ et } \quad 
\rho_X \big|_{\Gamma_Y} = \rho_Y \enspace.
$$
En particulier, $X\subset Y $ entra"ne $T(\Gamma_X) \subset T(\Gamma_Y)$. 
De plus, $Y$ est la variŽtŽ de torsion minimale contenant $X$ si et seulement si $\rho_Y(\Gamma_Y)\subset \mug^\infty$ et  $\rho_X(\lambda) \notin \mug^\infty$ pour
tout $\lambda \in \Gamma_X\setminus \Gamma_Y$. 
Autrement-dit, cette variŽtŽ de torsion minimale est caractŽrisŽe par   
$$
\Gamma_Y = \{\lambda\in\Gamma_X : \rho_X(\lambda) \in\mug^\infty\} \quad \mbox{ et }
\quad \rho_Y= \rho_X\big|_{\Gamma_Y} \enspace. 
$$ 

Le degrŽ d'obtruction des translatŽs de sous-tores dans des sous-tores
relativement ˆ $\k$ se rŽalise par un diviseur qui lui aussi est un translatŽ de sous-tore: 
 
\begin{prop}\label{cosettor}
Soient $X, Y \subset \G_m^N$ des translatŽs de sous-tores tels que $X \varsubsetneq Y$, alors il existe $\lambda\in\Gamma_X \setminus\Gamma_Y$ tel que 
$$
\omega_\k( X; Y)= \deg(\div(x^\lambda-\rho_X(\lambda))\cdot Y) \enspace. 
$$
\end{prop}

\begin{demo}
Soient $I(X), I(Y) \subset \k[x_1^{\pm1},\dots,x_N^{\pm1}]$ les idŽaux de dŽfinition de $X$ et de $Y$ respectivement. L'anneau $\k[Y]$ est factoriel car isomorphe {\`a} un anneau du type $\k[y_1^{\pm1},\dots,y_p^{\pm1}]$, et donc tout diviseur est principal. 
Soit donc $F\in I(X) \setminus I(Y)$ rŽalisant $\omega_{\k}(X;Y)$; on peut Žcrire $F$ modulo $I(Y)$ sous la forme
$$F = \sum_{\nu\in\Z^N/\Gamma_Y} f_\nu x^\nu 
= \sum_{\sigma\in\Z^N/\Gamma_X} x^\sigma F_\sigma
$$
avec $F_\sigma = \sum_{\tau \in \Gamma_X/\Gamma_Y}f_{\sigma+\tau} x^{\tau} $ pour 
$\sigma\in{\bf Z}^N/\Gamma_X$. Les mon™mes $x^\sigma$ ($\sigma\in\Z^N/\Gamma_X$) sont linŽairement indŽpendants modulo $I(X)$ et donc, comme $F\in I(X)$, on a $F_\sigma \in I(X)$ pour tout $\sigma$, de plus il existe $\sigma_0 \in \Z^N/\Gamma_X$ tel que 
$ F_{\sigma_0} \ne 0$. 
Ce polyn™me ne peut pas tre rŽduit ˆ un mon™me car $I(X)$ ne contient pas de mon™me et il existe donc $a,b\in \Gamma_X/\Gamma_Y$, $a\ne b$, tels que 
$$
\Newton(F) \supset \ov{ab} \enspace, 
$$
o $\ov{ab} \subset \R^N$ dŽsigne le segment d'extrŽmitŽs $a,b$. Soit $\lambda:=a-b$, on dŽduit de la proposition~\ref{degdiv} et de la monotonie du volume mixte par rapport ˆ l'inclusion~\cite[Thm.~4.12]{Ewa96}, que le degrŽ de $\div(F_{\sigma_0})\cdot Y$ est supŽrieur ou Žgal au degrŽ du diviseur $\div(x^\lambda-\rho_X(\lambda))\cdot Y$ et donc que $\omega_{\k}(X;Y)$ est rŽalisŽ par un diviseur qui est un translatŽ de sous-tore de la forme indiquŽe.
\end{demo}

On notera que pour un sous-corps $K$ qui n'est pas algŽbriquement clos, le degrŽ d'obstruc\-tion d'un translatŽ de sous-tore dans un sous-tore n'est pas nŽcessairement rŽalisŽ par un diviseur qui est un translatŽ de sous-tore. Voici un exemple, aimablement fourni par David~: 

\begin{exmpl}
Soit $\xi \in \Qbar\setminus \Q$ et posons 
$$
X:=(\xi,1-\xi )\in\G_m^2 \enspace \mbox{ et } \enspace Y:=\G_m^2 \enspace, 
$$ 
On a $\omega_{\bf Q}(X;Y)=1$ rŽalisŽ par la droite d'Žquation $x+y=1$, et on vŽrifie qu'il n'y pas d'Žquation binomiale de degrŽ $1$, dŽfinie sur $\Q$ et contenant $X$. Par ailleurs, on a $\omega_{\overline{\bf Q}}(X;Y)=1$ qui est bien sžr rŽalisŽ 
par la droite d'Žquation $x+y=1$, mais aussi par des diviseurs translatŽs de 
sous-tores, d'Žquations $x-\xi=0$,  $y-(1-\xi)=0$ ou encore $(1-\xi)x-\xi y=0$.
\end{exmpl}

Soient $X\subset Y$ des translatŽs de sous-tores de dimension $n$ et $p$ respectivement 
et posons 
$$
\Gamma_{X,Y} := \Gamma_X/\Gamma_Y
$$ 
le module quotient, qui est un $\Z$-module libre de rang $(N-n)-(N-p)=p-n$. ConsidŽrons aussi l'espace lin{\'e}aire $\Gamma_{X,Y}^\R:=\Gamma_{X,Y}\otimes \R$ muni de la m{\'e}trique quotient $||\cdot||_\bot$ induite de la mŽtrique euclidienne $\Vert\cdot\Vert_2$ par identification avec l'orthogonal de $\Gamma_Y^\R$ dans $\Gamma_X^\R$. En particulier $\Gamma_{X,Y}$ est un r{\'e}seau de $\Gamma_{X,Y}^{\bf R}$.

\begin{prop} \label{corbeau}
Soient $X\subsetneq Y$ des translatŽs de sous-tores de dimension $n$ et $p$ respectivement, alors pour tout $\lambda\in\Gamma_X \setminus \Gamma_Y$ on a
$$
\left(\binom{N}{p}\binom{N}{p-1}\right)^{-1/2} \leq \frac{\deg\Big(
\div(x^\lambda-\rho_X(\lambda))\cdot Y) \Big)}{||\lambda||_\bot\deg(Y)} 
\leq \left(\binom{N}{p}\binom{N}{p-1}\right)^{1/2} \enspace.
$$
\end{prop}

\begin{demo}
Posons $W:=\div(x^\lambda-\rho_X(\lambda))\cdot Y$, quitte \`a diviser $\lambda$ par un entier, ce qui divise \'egalement $\deg(W)$ et $\Vert\lambda\Vert_\bot$ par ce m\^eme entier, on peut supposer sans perte de gŽnŽralitŽ que $\Gamma_W:=\Gamma_Y+\Z\, \lambda$ est un $\Z$-module satur{\'e} et donc~\cite[\S~3]{PS05}
$$
\binom{N}{p-1}^{-1/2} \Vol_{N-p+1}(\Gamma_W^\R/\Gamma_W) \leq \deg(W) \leq \binom{N}{p-1}^{1/2} \Vol_{N-p+1}(\Gamma_W^\R/\Gamma_W)\enspace.
$$
De m\^eme
\begin{equation}\label{encadrementdegvol}
\binom{N}{p}^{-1/2} \Vol_{N-p}(\Gamma_Y^\R/\Gamma_Y) \leq \deg(Y) \leq \binom{N}{p}^{1/2} \Vol_{N-p}(\Gamma_Y^\R/\Gamma_Y)\enspace.
\end{equation}
Mais, $\Vol_{N-p+1}(\Gamma_W^\R/\Gamma_W)=||\lambda||_\bot\,
\Vol_{N-p}(\Gamma_Y^\R/\Gamma_Y)$ et les in\'egalit\'es cherch\'ees en r\'esul\-tent.
\end{demo}

\begin{lem} \label{vautour}
Avec les notations de la proposition~\ref{corbeau}, le premier minimum du rŽseau $\Gamma_{X,Y}$ est majorŽ par
$$\frac{2}{\sqrt{\pi}}\left(\binom{N}{n}\binom{N}{p}
\Gammag\left(1+\frac{p-n}{2}\right)^2 \right)^{\frac{1}{2(p-n)}} \left(\frac{\deg(X)}{\deg(Y)}\right)^\frac{1}{p-n}\enspace;$$
o $\Gammag$ dŽsigne la fonction gamma d'Euler. 
\end{lem}

\begin{demo}
Le th{\'e}or{\`e}me de Minkowski~\cite[Thm.V, \S~VIII.4.3]{Cas71} entra\^\i ne que le premier minimum de la norme $\Vert\cdot\Vert_\bot$ sur le rŽseau $\Gamma_{X,Y}$ est majorŽ par
$$
\frac{2}{\sqrt{\pi}} \Gammag\left(1+\frac{p-n}{2}\right)^{\frac{1}{p-n}}\, \Vol_{p-n}(\Gamma_{X,Y}^\R/ \Gamma_{X,Y})^\frac{1}{p-n}
\enspace.$$
Mais d'aprs~\cite[\S~3]{PS05}, {\it voir} inŽgalitŽ~(\ref{encadrementdegvol}) ci-dessus,
$$
\Vol_{p-n}(\Gamma_{X,Y}^\R/ \Gamma_{X,Y}) = 
\frac{\Vol_{N-n}(\Gamma_{X}^\R/\Gamma_{X})}{\Vol_{N-p}(\Gamma_{Y}^\R/\Gamma_{Y})}
\leq \left(\binom{N}{n}\binom{N}{p}\right)^{\frac{1}{2}} \frac{\deg(X)}{\deg(Y)}
\enspace, 
$$
d'o l'{\'e}nonc{\'e}.
\end{demo}

En combinant le lemme~\ref{vautour} avec l'inŽgalitŽ de droite de la proposition~\ref{corbeau} on trouve:

\begin{cor} \label{omega}
Soient $X\subsetneq Y$ des translatŽs de sous-tores de dimension $n$ et $p$ respectivement, alors
$$
\omega_\k(X;Y) \leq \use{ctedsomega}(N,p,n) \deg(X)^{\frac{1}{p-n}} \deg(Y)^{1-\frac{1}{p-n}}
$$
avec
$\declare{ctedsomega}(N,p,n) := \frac{2}{\sqrt{\pi}} \left(\binom{N}{p}\binom{N}{p-1}\right)^{\frac{1}{2}} \left(\binom{N}{n}\binom{N}{p}\Gammag\left(1+\frac{p-n}{2}\right)^2 \right)^{\frac{1}{2(p-n)}} \leq 4^{N}{N}^{3/2}$.
\end{cor}

En particulier lorsque $X$ est rŽduit ˆ un point $\alpha\in \G_m^N$ on a  $\omega_\k(\alpha;Y) \leq (2^NN)^{3/2}\deg(Y)^{1-\frac{1}{p}}$.

\medskip
Il est naturel de se demander si cette majoration s'Žtend ˆ d'autres vari\'et\'es, \cad  si  
$$
\omega_\k(X,Y) \le c(N) \deg(X)^\frac{1}{p-n} \, \deg(Y)^{1-\frac{1}{p-n}}
$$
pour des variŽtŽs projectives $X\varsubsetneq Y $ de dimension $n$ et $p$ respectivement. M.~Chardin~\cite{Cha89} a d{\'e}montr{\'e} un tel r\'esultat pour $Y=\P^N$ et $X\subset{\bf P}^N$ quelconque. Cependant, il semblerait qu'une gŽnŽralisation ne soit possible qu'en imposant des conditions sur la r\'egularit\'e de Castelnuovo-Mumford de $Y$.

\medskip
Il serait Žgalement int{\'e}ressant de remplacer les estimations de la proposition~\ref{corbeau} par une {\'e}galit{\'e}. Est-ce possible en modifiant convenablement la m{\'e}trique de $\Gamma_{X,Y}$? Le r{\'e}sultat suivant fourni une expression alternative {\it exacte} pour le degrŽ d'obstruction, cependant elle ne nous est pas utile pour la suite.

\begin{prop} \label{carpincho}
Soit $Y\subset \G_m^N$ un translatŽ de sous-tore, $\cB \in (\Z^p)^N$ satisfaisant $L_\cB=\Z^p$ et $\beta\in \G_m^N$, tels que $Y=\beta T_\cB$. Soit $\lambda\in\Z^N\setminus\Gamma_Y$ tel que $\lambda \Z+\Gamma_Y$ soit satur{\'e} et $\gamma\in \k^\times$ quelconque, alors
$$
\deg( \div(x^\lambda -\gamma)\cdot Y 
)= (p-1)! \cdot \Vol_{p-1}(\pi_\lambda(Q_\cB))\cdot ||M_\cB(\lambda)||_2
\enspace,
$$
o{\`u} 
$Q_\cB \subset \R^p $ d{\'e}signe l'enveloppe convexe des points ${\mathbf 0},b_1,\dots,b_N$, $\pi_\lambda$ la projection orthogonale $\pi_\lambda:\R^p\to M_\cB(\lambda)^\bot \cong \R^{p-1}$ et $||\cdot||_2$ la norme euclidienne.
\end{prop}


\begin{demo}
Posons $W:=\div(x^\lambda -\gamma)\cdot Y$, comme $\lambda \Z+\Gamma_Y$ est satur{\'e} $W$ est rŽduit et irrŽductible de dimension $p-1$: c'est 
donc un translatŽ de sous-tore. 

\smallskip 
Soit $E:= M_\cB(\lambda)^\bot$ et  $E^\Z:=  \pi_\lambda(\Z^p)$, alors $E^\Z$ est un rŽseau de $E$ et on vŽrifie ${\Vol_{p-1}(E/E^\Z)} = ||M_\cB(\lambda)||_2^{-1}$. Posons maintenant $c_j:=\pi_\lambda(b_j)$  pour $j=1,\dots, N $ et $\cC:=(c_1,\dots,c_N)$, de sorte que $L_\cC=E^\Z$, et soit $\eta\in \G_m^N$ tel que $\gamma=\eta^\lambda$, de sorte que $\eta T_\cC$ est de dimension $p-1$. Le translatŽ de sous-tore $\eta T_\cC$ est l'image de l'application $s\mapsto (\eta_1\,s^{c_1},\dots,\eta_N \,s^{c_N})$ et on vŽrifie 
$$
(\eta_1\,s^{c_1},\dots,\eta_N \,s^{c_N})^\lambda -\eta^\lambda
=\eta^\lambda\, ( s^{c_1\, \lambda_1+\cdots+c_N\,\lambda_N}-1)\enspace.
$$
Mais $c_1\, \lambda_1+\cdots+c_N\,\lambda_N=\pi_\lambda(b_1)\, \lambda_1+\cdots +\pi_\lambda(b_N)\, \lambda_N =\pi_\lambda(M_\cB(\lambda))=0 $ d'o{\`u} $\eta T_{\cC}\subset W$ et donc $\eta T_{\cC}=W$. Le calcul de degr{\'e} dŽcoule alors du calcul du volume du polytope $P:= \Conv({\mathbf 0},c_1,\dots,c_N) =\pi_\lambda(Q_\cB)$:
$$
\deg(W)= (p-1)!\frac{\Vol_{p-1}(P)}{\Vol_{p-1}(E/E^\Z)}
= (p-1)! \cdot \Vol_{p-1}(\pi_\lambda(Q_\cB))\cdot ||M_\cB(\lambda)||_2
\enspace.
$$
\end{demo}

\begin{rem}
En combinant les propositions~\ref{carpincho} et~\ref{degdiv} on vŽrifie
$$\MV\big([0,M_\cB(\lambda)], \underbrace{Q_\cB,\dots,Q_\cB}_{p-1\ fois}\big) = (p-1)!{\Vol}_{p-1}\left(\pi_\lambda(Q_\cB)\right)\cdot\Vert M_\cB(\lambda)\Vert_2
\enspace.$$
Cette identitŽ se dŽduit Žgalement de la formule gŽnŽrale suivante pour un convexe $Q$ quelconque en lieu et place du segment $[0,M_\cB(\lambda)]$~:
$$\MV\big(Q, \underbrace{Q_\cB,\dots,Q_\cB}_{p-1\ fois}\big) = (p-1)!\sum_{w\in S^{p-1}}a_w(Q)\Vol_{p-1}(Q_w)
\enspace,$$
o $w$ parcourt la sphre unitŽ de $\R^p$, $a_w(Q):=\max\{\langle w,q\rangle:q\in Q\}$ et $Q_w:=\{q\in Q:\langle w,q\rangle=a_w(q)\}$, {\it voir}~\cite[Ch.~IV, Thm.~4.10]{Ewa96}.
\end{rem}

\section{D{\'e}monstrations du th{\'e}or{\`e}me~\ref{laurier} et des corollaires~\ref{corCM} et~\ref{bertrand}}\label{demonstrations}

Nous donnons ici les preuves des rŽsultats ŽnoncŽs dans l'introduction. Remarquons tout de mme que dans le cas des hypersurfaces ($n=N-1$) le thŽorme~\ref{laurier} est trivial: dans ce cas un translatŽ de sous-tore qui n'est pas de torsion est de la forme $X={Z}(x^b-\lambda)$ pour un vecteur $b\in \Z^N$ ˆ coordonnŽes premires entre elles et 
$\lambda\in K^\times\setminus \mug^\infty$, et on a $U_X=\G_m^N$. Donc, d'aprs~\cite[prop.~VI.4]{PS05}
$$
\wh{\mu}^\ess(X) = \frac{\hnorm(X)}{\deg(X)} = \frac{\hnorm(\lambda)}{\omega_{\overline{\Q}}(X;\G_m^N)} \geq  m(K) \frac{\deg(U_X)}{\omega_{\overline{\Q}}(X;U_X)}\enspace,
$$
o $\hnorm(X)$ dŽsigne la hauteur normalisŽe de l'hypersurface $X$.

\begin{demo}[D{\'e}monstration du th{\'e}or{\`e}me~\ref{laurier}]
On Žcrit $X=\alpha T$ et $U=U_X=\xi T'$ pour certains sous-tores $T$, $T'$, $\alpha\in (K^\times)^N$ et $\xi\in(\mug^\infty)^N$. On pose $n=\dim(X)$, $p=\dim(U_X)$ et on suppose sans perte de gŽnŽralitŽ $n<p$.
\smallskip

Soit $\lambda \in\Gamma_X\setminus \Gamma_U$, \'ecrivons $\lambda= c+b $ avec $b \in \Gamma_U^\Q$ et $c \in(\Gamma_U^\Q)^\bot$. Soit $\zeta=(\zeta_1,\dots,\zeta_N)\in X$ un point quelconque, alors $\zeta=\alpha\theta$ pour un certain $\theta\in T_X$, d'o, puisque $\theta^\lambda = 1$,
$$
\zeta^\lambda= (\alpha\theta)^\lambda = \alpha^\lambda \quad  \mbox{ pour tout } \enspace \lambda \in \Gamma_X 
\enspace.
$$
De plus $\zeta^b\in\mug^\infty$ car $\zeta\in\xi T_U$, $b\in \Gamma_U^\Q$ et pour tout entier $m$ tel que $mb\in\Gamma_U$ on a $(\zeta^{b})^m=\xi^{mb}$. Comme $\lambda\notin\Gamma_U$ il vient que $\zeta^\lambda$ (et par la suite $\zeta^c$) n'est pas une racine de l'unitŽ, {\`a} cause de la minimalitŽ de $U$ parmi les variŽtŽs de torsion contenant $X$.
\smallskip

Soit $L \supset K$ un corps de nombres contenant $\zeta_1^{c_1}, \dots, \zeta_N^{c_N}$. Comme $\alpha^\lambda$ et $\zeta^c$ ne diffrent que du facteur $\zeta^b$ qui est une racine de l'unitŽ, pour toute place $v \in M_L$ on a 
$$\begin{array}{rcl}
\big|\log|\alpha^\lambda|_v\big| &= &\big|\log|\zeta^c|_v\big|\\[2.5mm]
&= &\Big| c_1\,\log|\zeta_1|_v +\cdots
+ c_N\,\log|\zeta_N|_v\Big|\\[2.5mm]
&\leq &||c||_1\, \max\Big(\big|\log|\zeta_1|_v\big|, \dots, \big|\log|\zeta_N|_v\big| \Big)\\[2.5mm]
&\leq &\sqrt{p-n}\,||\lambda||_\bot \, \max\Big(\big|\log|\zeta_1|_v\big|, \dots, \big|\log|\zeta_N|_v\big| \Big)
\enspace,\end{array}
$$
car $c$ est la projection de $\lambda$ sur l'orthogonal $(\Gamma_U^{\bf Q})^\bot$ de $\Gamma_U^{\bf Q}$ dans $\Gamma_X^\Q$ et $\dim_{\bf Q}(\Gamma_U^{\bf Q})^\bot = p-n$. D'o\`u, en sommant sur toutes les places de $L$
$$
2\,\wh{h}(\alpha^\lambda) = \sum_{v\in M_L} \frac{[L_v:\Q_v]}{[L:\Q]}\, 
\big|\log|\alpha^\lambda|_v\big| \leq \sqrt{p-n}\,||\lambda||_\bot \,(\hnorm(\zeta) + \hnorm(\zeta^{-1}))
\enspace.$$
Et comme 
$$
\wh{h}(\zeta^{-1}) = \sum_v \frac{[L_v:\Q_v]}{[L:\Q]} \max\{ 0, -\log|\zeta_1|_v, \dots, -\log|\zeta_N|_v \} \le \wh{h}(\zeta_1^{-1})+\cdots + \wh{h}(\zeta_N^{-1}) \leq N \wh{h}(\zeta) \enspace, 
$$ 
\begin{equation}\label{minohnorm}
\hnorm(\zeta) \geq \frac{2}{(N+1)\sqrt{p-n}}\frac{\wh{h}(\alpha^\lambda)}{\Vert\lambda\Vert_\bot} \geq \frac{2}{(N+1)\sqrt{p-n}}\frac{ m(K)}{\Vert\lambda\Vert_\bot}
\end{equation}
par dŽfinition de $m(K)$, car $\alpha^\lambda\in K^\times\setminus\mug^\infty$. 
Soit maintenant $\lambda\in\Gamma_X\setminus\Gamma_U$ r{\'e}alisant le degrŽ
d'obstruction $\omega_\Qbar(X;U)$, par la proposition~\ref{corbeau} 
$$
||\lambda||_\bot \leq \left(\binom{N}{p}\binom{N}{p-1}\right)^{1/2}\frac{\omega_\Qbar (X;U)}{\deg(U)}
 \enspace,
$$
d'o{\`u} finalement
$$
\hnorm(\zeta) \geq c(N,p,n) \cdot 
 m(K) \cdot  \frac{\deg(U)}{\omega_{\overline{\Q}} (X;U)}
$$
avec $c(N,p,n) = 2\left((N+1)^2(p-n)\binom{N}{p}\binom{N}{p-1}\right)^{-1/2} \ge N^{-3/2} 2^{-N}$, ce qui achve la dŽmonstration 
du thŽorme~\ref{laurier}.
\end{demo}

\begin{rem}\label{remsinnou}
En minorant $ m(K)$ par $c[K:{\bf Q}]^{-1}$ selon la rŽponse attendue ˆ la question de Lehmer, on est tentŽ de considŽrer le degrŽ d'obstruction suivant, pour $X$ et $U$ comme
dans la dŽmonstration ci-dessus: 
$$
\inf\left\{[{\bf Q}(\alpha^\lambda):{\bf Q}]
\deg(\div(x^\lambda-\alpha^\lambda)\cdot U) : \enspace \lambda\in\Gamma_{X}\setminus\Gamma_{U}\right\}\enspace.
$$
\end{rem}

\begin{lem}\label{translatecasCM}
Soit $X$ un translatŽ de sous-tore dŽfini sur un corps CM, alors la plus petite variŽtŽ de torsion $G_X$ contenant $X\cdot \ov X$ est un sous-tore. 
\end{lem}

\begin{demo}
Soit $K$ le corps de dŽfinition de $X$, Žcrivons $X=\alpha T$ et $G_X=\xi T'$ pour certains sous-tores $T$ et $T'$, $\alpha \in (K^\times)^N$ et  $\xi\in(\mug^\infty)^N$. On a $X\cdot\overline{X}=|\alpha|^2 T$ donc pour tout $\lambda\in\Gamma_{T'}$ on a 
$|\alpha|^{2\lambda} = \xi^\lambda\in {\bf R}_+ \cap \mug^\infty$ d'o $\xi^\lambda=1$; autrement-dit $\xi\in T'$ et donc $G_X = \xi T' = T'$ est un sous-tore.
\end{demo}

\begin{demo}[DŽmonstration du corollaire~\ref{corCM}]
Lorsque $K$ est un corps totalement rŽel ou une extension abelienne de ${\bf Q}$ on a $m(K)\geq \frac{1}{2}\log\left(\frac{1+\sqrt{5}}{2}\right)$ d'aprs~(\ref{schinzel}) ou $m(K)\geq \frac{\log(5)}{12}$ d'aprs~(\ref{amodvo}) et on reporte directement dans le thŽorme~\ref{laurier}.

Si $K$ est un corps CM, la plus petite sous-variŽtŽ de torsion contenant $Y:=X\cdot\overline{X}$ est un sous-tore $G_X$, d'aprs le lemme~\ref{translatecasCM}. De plus, $Y$ est dŽfinie sur le sous-corps totalement rŽel $K_0$ de $K$. Mais on a $m(K_0)\geq \frac{1}{2}\log\left(\frac{1+\sqrt{5}}{2}\right)$ d'aprs~(\ref{schinzel}) et, si $Y$ n'est pas une variŽtŽ de torsion, on dŽduit du thŽorme~\ref{laurier}~:
$$\munorm^\ess(Y) \geq \frac{\use{ctedslaurier}(N)}{2}\log\left(\frac{1+\sqrt{5}}{2}\right) \frac{\deg(G_X)}{\omega_{\overline{\bf Q}}(Y;G_X)}
\enspace.$$

Maintenant, si $Y$ est une variŽtŽ de torsion on a $Y=G_X$ et $\omega_{\overline{\bf Q}}(Y;G_X)=\infty$, la minoration ci-dessus est donc encore vraie dans ce cas. Finalement on remarque $\munorm^\ess(Y) \leq 2\munorm^\ess(X)$ pour conclure la dŽmonstration.
\end{demo}


\begin{demo}[DŽmonstration du corollaire~\ref{bertrand}]
On combine~(\ref{minohnorm}) dans la preuve du thŽorme~\ref{laurier} ci-dessus et le lemme~\ref{vautour}, ce qui donne ($X=\alpha T$, $U_X=\xi T'$)~:
\begin{eqnarray*}
\munorm^{\abs}(\alpha T) &\geq &\frac{2}{(N+1)\sqrt{p-n}}\cdot  m(K)\cdot\frac{1}{\Vert\lambda\Vert_\bot}\\[3mm]
&\geq &\frac{\sqrt{\pi}}{(N+1)\sqrt{p-n}}\left(\binom{N}{n}\binom{N}{p} \Gammag\left(1+ \frac{p-n}{2}\right)^{2}\right)^{-\frac{1}{2(p-n)}} \cdot  m(K) \cdot \left(\frac{\deg(T')}{\deg(T)}\right)^{\frac{1}{p-n}}\\[3mm]
&\geq &\declare{ctedsbertrand}(N)\cdot m(K)\cdot \left(\frac{\deg(U_X)}{\deg(X)}\right)^{\frac{1}{p-n}}
\enspace,
\end{eqnarray*}
avec $\use{ctedsbertrand}(N)\geq 2^{-N}N^{-2}$.
\end{demo}

\begin{rem}\label{diffomegas}
Notons que la minoration de la conjecture~\ref{sinnoubogo}, si elle est de nature purement gŽomŽtrique, n'en est pas toujours meilleure que celle de la conjecture~\ref{sinnou}. Par exemple, dans ${\bf G}_m^3$, soit $X$ une courbe arbitraire dans le plan d'Žquation $x=a$ avec $a\in \Q^\times \setminus\{\pm1\}$. Alors $V_X = Z(x-a)$ et  $U_X^\Q = \G_m^3$, d'o $\deg(V_X)=\deg(U_X^\Q)=1$ puis $\omega_{\bf Q}(X;U^\Q_X) = 1$ tandis que $\omega_{\overline{\bf Q}}(X;V_X) = \deg(X)$ est arbitraire et donc 
$$
\frac{\deg(U_X^\Q)}{\omega_{\bf Q}(X;U_X^\Q)} \gg 
\frac{\deg(V_X)}{\omega_{\overline{\bf Q}}(X;V_X)} \enspace.
$$ 
\end{rem}

\section{\smash{\`A} propos du problme de Lehmer sur ${\bf k}(t)$}\label{lehmerfonct}

Soit ${\bf k}$ un corps algŽbriquement clos et $t$ une variable. Pour une extension finie $K$ de ${\bf k}(t)$ on considre un revtement de la droite projective 
$$
\pi: C_K\to \P^1(\k)
$$
de degrŽ $\deg(\pi)=[K:\k(t)]$ par une courbe projective lisse $C=C_K$, correspondant ˆ l'extension de corps $\k(t)\hookrightarrow K$. Un ŽlŽment $\xi \in K$
s'interprte ainsi comme une fonction rationnelle $\xi:C\dashrightarrow \k$, et
pour une place $v\in C$ on pose $\Ord_v(\xi)\in \Z$ l'ordre d'annulation
de $\xi$ en $v$. On a la \og formule du produit\fg, qui s'Žcrit additivement  
$$
\sum_{v\in C}\Ord_v(\xi)=0 \enspace,
$$
et la {\em hauteur}  d'un point  $\alpha=(\alpha_0:\cdots:\alpha_N)\in {\bf P}^N(K)$ est par dŽfinition
$$
h(\alpha) := \frac{1}{[K:{\bf k}(t)]} \sum_{v\in C}
\max\Big(-\Ord_v(\alpha_0),\dots,-\Ord_v(\alpha_N)\Big)
\enspace.$$
Cette formule est indŽpendante du choix des coordonnŽes projectives de $\alpha$ 
gr‰ce ˆ la formule du produit, et elle ne dŽpend pas non plus du choix du revtement $\pi$.

Pour  $\alpha \in \P^N(k(t))$ dont les coordonnŽes $\alpha_j $ sont des polyn™mes premiers entre eux, on a $h(\alpha)=\max_j \deg(\alpha_j)$. En outre, $h(\alpha)=0$ si et seulement si $\alpha$ est \og constant\fg\ par rapport ˆ $t$, \cad si et seulement si 
$\alpha\in {\bf P}^N({\bf k})$: on le vŽrifie en supposant  par exemple $\alpha_0\ne 0$, alors $h(\alpha)=0$ si et seulement si $\Ord_v(\alpha_j/\alpha_0)\geq 0$ pour $j=1,\dots,N$ et tout $v\in C$, ce qui Žquivaut ˆ ce que $\alpha_j/\alpha_0\in \k$. Ainsi, les ŽlŽments de $\k$ jouent dans le cas fonctionnel, le r™le des racines de l'unitŽ dans le cas arithmŽtique.

\medskip
L'application $\alpha:C \dashrightarrow \P^N(\k)$ s'Žtend en une section globale $C\to\P^N(\k)$ car $C$ est supposŽe lisse. Notons $\cX_\alpha\subset C\times \P^N(\k)$ son graphe; on a~\cite[\S~B.10]{HS00}
$$
h(\alpha)= \frac{1}{[K:\k(t)]} \Card\Big( \cX_\alpha \cap
p_2^{-1}(E)\Big) \enspace, 
$$
o $p_2$ dŽsigne la projection $  C\times \P^N(\k) \to \P^N(\k)$ et $E$ est un hyperplan gŽnŽrique de  $\P^N(\k)$.

\medskip
Maintenant considŽrons le cas d'une variŽtŽ  $X\subset{\bf P}^N(\overline{{\bf k}(t)})$ de dimension quelconque $d\ge 0$, dŽfinie sur $K$. ConsidŽrons la variŽtŽ $\cX \subset C\times\P^N(\k)$ de dimension $d+1$ dont la fibre gŽnŽrique sur $C$ est Žgale ˆ $X$. Soit $C\hookrightarrow \P^M(\k)$ une immersion de $C$ dans un espace projectif. Soient 
$$
p_1:C\times \P^N(\k)\to C \subset \P^M(\k) \quad ,\quad p_2: C\times
\P^N(\k) \to \P^N(\k)
$$
les deux projections naturelles,  
et pour $d_1,d_2\ge 0$ tels que $d_1+d_2=d+1$ considŽrons le
multidegrŽ
$$
\deg_{(d_1,d_2)}(\cX) := \Card(\cX \cap p_1^{-1}(E_1) \cap
p_2^{-1}(E_2))
$$
o $E_1 \subset \P^M$ et $E_2\subset \P^N$ dŽsignent des espaces linŽaires 
gŽnŽriques de codimension $d_1$ et $d_2$, respectivement.
La {\it hauteur} de la variŽtŽ $X$ est alors par dŽfinition 
$$
h(X) := \frac{1}{[K:\k(t)]}\deg_{(0,d+1)}(\cX) \enspace.
$$
On peut vŽrifier que ceci ne dŽpend pas des choix effectuŽs, et qu'elle co•ncide avec la dŽfinition prŽcŽdente dans le cas des points. On a $h(X)\ge 0$ et de fait $h(X)=0$ si et seulement si $p_2(\cX)$ est de dimension $d$, ce qui Žquivaut ˆ ce que $X$ soit dŽfinie sur $\k$. Pour les autres multidegrŽs, on a 
$$
\deg_{(d_1,d_2)}(\cX) = \left\{ \begin{array}{ll} 
[K:\k(t)]\deg(X) & \mbox{ pour } d_1=1, d_2=d \enspace, \\[2mm] 
 0 & \mbox{ pour } d_1\ge 2 \enspace.
\end{array} \right.
$$

La hauteur et le degrŽ de $X$ peuvent s'interprŽter en termes de formes 
rŽsultantes de l'idŽal bihomogne
$I(\cX)\subset \k[w_0,\dots,w_M,x_0,\dots,x_N]$. 
Soient $W$ un groupe de $M+1$ variables 
et $V_0,\dots,V_d$ des groupes de $N+1$ variables chacun et soit
$$
\res_\cX \in \k[W ,V_0,\dots,V_d]
$$
la forme rŽsultante de $I(\cX)$ d'indice $(1,0)$ et $(0,1),\dots,(0,1)$ 
($d+1$ fois), 
{\it voir}~\cite[\S~3]{Rem01a} ou encore~\cite[\S~I.2]{PS04} 
pour la d{\'e}finition et propri{\'e}t{\'e}s de base des formes r{\'e}sultantes. 
C'est un polyn™me multihomogne dont le degrŽ par rapport ˆ chaque groupe 
de variables est~\cite[Prop.~3.4 et 2.11]{Rem01a}
$$\begin{array}{rclcl}
\deg_{W}(\res_\cX) &= &\deg_{(0,d+1)}(\cX) &= &[K:\k(t)] h(X) \enspace,\\[2mm]
\deg_{V_i}(\res_\cX)&= &\deg_{(1,d)}(\cX) &= &[K:\k(t)] \deg(X)  \enspace, \quad 
\mbox{ pour } i=1,\dots, d \enspace. 
\end{array}$$
Comme consŽquence de cette interprŽtation et de~\cite[Lem.~2.11]{Rem01a}
on peut dŽmontrer le thŽorme de BŽzout fonctionnel: 
soit $f\in \k[w_0,\dots,w_M,x_0,\dots,x_N] \setminus I(\cX)$ 
une forme bihomogne,  alors
\begin{equation} \label{bezout}
h(X\cdot \div(f)) =  \deg_\x (f) h(X)+ \deg_\w(f)\deg(X)\enspace.
\end{equation}
Similairement, on obtient une formule de Hilbert-Samuel pour la dimension de l'espace des formes de bi-degrŽ $(\eta,\delta)$ modulo $I(\cX)$~\cite[Thm.~2.10]{Rem01a}: 
\begin{eqnarray} \label{hilbert}
\dim_\k \Big( \k[\cX]_{(\eta,\delta)} \Big)
&=&\frac{\deg_{(0,d+1)}(\cX)}{(d+1)!} \delta^{d+1} + \frac{\deg_{(1,d)}(\cX)}{d!}\eta\delta^d + O((\eta+\delta)\delta^{d-1}) \\[2mm]
&=&[K:\k(t)]\left(\frac{h(X)}{(d+1)!} \delta^{d+1} + \frac{\deg(X)}{d!}\eta\delta^d\right)(1+o(1))
\nonumber 
\end{eqnarray}
pour $\eta,\delta \to \infty$.

Par ailleurs, le {\it minimum essentiel} de $X$ est naturellement dŽfini par 
$$
\mu^{\ess}(X) := \inf\Big\{ \theta : \{ \xi\in X : h(\xi) \le \theta\} \mbox{ est Zariski dense}\Big\}
\enspace.
$$
On pose $U_X\subset \P^N(\ov{\k(t)})$ la plus petite variŽtŽ dŽfinie sur $\k$ et contenant $X$. 
Le rŽsultat principal de ce paragraphe est l'analogue fonctionnel suivant 
de la conjecture~\ref{sinnou}: 

\begin{prop} \label{sinnoufonctionnel}
Soit $X\subset{\bf P}^N(\ov{{\bf k}(t)})$ une variŽtŽ de dimension $d\ge 0$, alors 
$$(d+1)\mu^{\ess}(X) \geq \frac{h(X)}{\deg(X)} \geq \frac{\deg(U_X)}{\omega_{{\bf k}(t)}(X;U_X)}\enspace.$$
\end{prop}

\begin{demo} 
La premire inŽgalitŽ rŽsulte de la proposition~\ref{zhang} ci-dessous. Pour la seconde, on remarque $U_X=p_2(\cX) \subset {\bf P}^N({\bf k})$. Cette variŽtŽ est donc 
Žgale ˆ $X$ si $X$ est dŽfinie sur $\k$, et de dimension $d+1$ sinon. 

Si $X=U_X$ on a $\omega_{{\bf k}(t)}(X;U_X)=\infty$ et le rŽsultat est clair, sinon $X$ est de codimension $1$ dans $U_X$. On a clairement
$$
\deg_{(0,d+1)}(\cX) = \deg(U_X)\deg(p_2\vert_{\cX}) \geq \deg(U_X)\enspace.$$
D'un autre c™tŽ, le diviseur minimal de $U_X$ dŽfini sur $\k(t)$ et contenant $X$ est 
$\bigcup_{\sigma} \sigma(X)$ o $\sigma$ parcours les $\k(t)$-immersions du corps de dŽfinition $K$ de $X$ dans la cl™ture algŽbrique $\ov{\k(t)}$. Le degrŽ de ce diviseur est $[K:\k(t)]\deg(X)$ et donc
$\omega_{{\bf k}(t)}(X;U_X)=[K:\k(t)]\deg(X)$. 
On en conclut 
$$
\frac{h(X)}{\deg(X)} = \frac{\deg_{(0,d+1)}(\cX)}{\deg_{(1,d)}(\cX)} \geq \frac{\deg(U_X)}{[K:\k(t)]\deg(X)} = \frac{\deg(U_X)}{\omega_{{\bf k}(t)}(X;U_X)}
\enspace.$$
\end{demo}

Il ne reste qu'ˆ dŽmontrer l'analogue fonctionnel suivant du thŽorme des minimums algŽbriques successifs. On reprend pour cela en partie les idŽes de~\cite[Thm.~3.1]{DaPh98}.

\begin{prop}\label{zhang}
Soit $X\subset{\bf P}^N(\overline{{\bf k}(t)})$ une variŽtŽ de dimension  $d\ge 0$, alors
$$\mu^{\ess}(X) \leq \frac{h(X)}{\deg(X)} \leq (d+1)\mu^{\ess}(X)\enspace.$$
\end{prop}

\begin{demo} Pour l'inŽgalitŽ de gauche, il suffit de dŽmontrer que pour tout diviseur $Z$ de $X$ il y a  des points dans $X\setminus Z$ de hauteur $\leq {h(X)}/{\deg(X)}$: en intersectant $X$ avec $d$ formes linŽaires gŽnŽriques ˆ coefficients dans ${\bf k}$ le  cycle intersection est de dimension $0$, de degrŽ $\deg(X)$ et ˆ support dans $X \setminus Z$. Par le thŽorme de BŽzout (formule~(\ref{bezout}) ci-dessus) la hauteur de ce cycle est $h(X)$, et donc son support contient au moins un point de hauteur $\leq {h(X)}/{\deg(X)}$. 

\medskip
Pour l'inŽgalitŽ de droite, on commence par supposer sans perte de gŽnŽralitŽ que $X$ est dŽfinie sur ${\bf k}(t)$. Soit $I(X)$ l'idŽal de dŽfinition de $X$ dans $\k(t)[x_0,\dots,x_N]$ et $\varepsilon>0$, il existe une base de $I(X)^\perp_\delta$ de la forme 
$$b_i=\big(\xi_i^\alpha\big)_{
\begin{array}{l}\scriptstyle\kern-1mm\alpha\in{\bf N}^{N+1}\\[-1.5mm]
\scriptstyle\kern-1mm|\alpha|=\delta\end{array}}
\enspace,\quad i=1,\dots,L:=\dim_{{\bf k}(t)}\big({\bf k}(t)[X]_\delta\big)\enspace,$$
avec $\xi_i\in X(\overline{{\bf k}(t)})$ et $h(\xi_i)\leq \mu^\ess(X)+\varepsilon$, par dŽfinition du minimum essentiel. On en dŽduit ˆ l'aide de la formule de Brill-Gordan~:
\begin{equation}\label{hauteurideal}
h(I(X)_\delta) = h(I(X)^\perp_\delta) \leq \sum_{i=1}^Lh(b_i) = \delta\sum_{i=0}^Lh(\xi_i) \leq \delta L(\mu^\ess(X)+\varepsilon)
\enspace.
\end{equation}

Par ailleurs on a
\begin{equation}\label{hilbertmorcele}
\begin{array}{rcl}
\dim_{\bf k}({\bf k}[\cX]_{(\eta,\delta)}) &= &\dim_{\bf k}\big({\bf k}[w,x]_{(\eta,\delta)}\big) - \dim_{\bf k}(I(\cX)_{(\eta,\delta)})\\[3mm]
&= &(\eta+1)\binom{\delta+n}{n} - \dim_{\bf k}(I(\cX)_{(\eta,\delta)})
\enspace.
\end{array}
\end{equation}
D'aprs~\cite[Cor.~2]{Thu95} ({\it voir} aussi~\cite[page~489]{Mah41}) il existe une base $c_1,\dots,c_M$ du ${\bf k}(t)$-espace $I(X)_\delta$ satisfaisant $\sum_{j=1}^Mh(c_j) = h(I(X)_\delta)$. Mais alors les ŽlŽments 
$$t^\ell c_j\enspace,\quad j=1,\dots,M\enspace,\quad \ell=0,\dots,\eta-h(c_j)$$
sont linŽairement indŽpendants sur ${\bf k}$ dans $I(\cX)_{(\eta,\delta)}$ qui est donc de dimension au moins $(\eta+1)M-\sum_{j=1}^Mh(c_j)\geq (\eta+1)M-h(I(X)_\delta)$. En reportant dans~(\ref{hilbertmorcele}) on obtient
\begin{equation}\label{hilbertrassemble}
\begin{array}{rcl}
\dim_{\bf k}({\bf k}[\cX]_{(\eta,\delta)}) &\leq &(\eta+1)\binom{\delta+n}{n}-(\eta+1)M+h(I(X)_{\delta})\\[3mm]
&\leq &(\eta+1)L+h(I(X)_{\delta})
\enspace.
\end{array}
\end{equation}
En rŽunissant~(\ref{hilbert}),~(\ref{hauteurideal}) et~(\ref{hilbertrassemble}) il vient avec l'expression asymptotique de $L$ comme fonction de Hilbert~:
$$
\left(\frac{h(X)}{(d+1)!}\delta^{d+1}+\frac{\deg(X)}{d!}\eta\delta^d\right)(1+o(1)) \leq \frac{\deg(X)}{d!}\delta^d\left(\eta+1+\delta(\mu^\ess(X)+\varepsilon)\right)(1+o(1))
$$
et enfin, en divisant par $\delta^{d+1}$ et en faisant tendre $\eta$ et $\delta$ vers l'infini de sorte que $\eta/\delta$ reste bornŽ,
$$
\frac{h(X)}{(d+1)!} \leq \frac{\deg(X)}{d!}\cdot(\mu^\ess(X) + \varepsilon)
\enspace,$$
puis le rŽsultat voulu lorsque $\varepsilon$ tend vers $0$.
\end{demo}


\let\livrefont=\sl

\typeout{References}

\end{document}